\documentclass[12pt]{amsart} 

\theoremstyle{plain}
\newtheorem{thm}{Theorem}[section]
\newtheorem*{theoremA}{Theorem A}
\newtheorem*{theoremB} {Theorem B}
\newtheorem{lem}[thm]{Lemma}
\newtheorem{prop}[thm]{Proposition}

\newtheorem{cor}[thm]{Corollary}
\theoremstyle{remark}
\newtheorem{rem}[thm]{Remark}
\theoremstyle{definition}
\newtheorem{definition}[thm]{Definition}

\DeclareMathOperator{\Co}{Co}
\DeclareMathOperator{\diam}{diam}
\begin{document}
\title[Approximately convex sets in normed spaces]{On the size of
approximately convex sets in normed spaces}
\author{S. J. Dilworth, Ralph Howard and James W. Roberts}
\address{Department of Mathematics, University of South Carolina,
 Columbia, SC 29208, U.S.A.} \email{dilworth@math.sc.edu,
howard@math.sc.edu, roberts@math.sc.edu} 
\thanks{The work of the second author was supported in part by DoD
Grant No. N00014-97-1-0806} 
\date{16 August 1999}
\begin{abstract} Let $X$ be a normed space. A set
$A \subseteq X$ is \textit{approximately convex} if
$d(ta+(1-t)b,A) \le 1$ for all $a,b \in A$
and $t \in [0,1]$.
 We prove that every $n$-dimensional normed space
contains approximately convex sets $A$ with $\mathcal{H}(A,\Co(A))\ge
\log_2n-1$ and $\diam(A) \le C\sqrt n(\ln n)^2$,
where $\mathcal{H}$ denotes the Hausdorff distance. These estimates are 
reasonably sharp. For every $D>0$, we construct 
worst possible approximately convex sets in 
$C(0,1)$ such that $\mathcal{H}(A,\Co(A))=\diam(A)=D$.
Several results pertaining to the Hyers-Ulam stability theorem are also
proved.
\end{abstract} 
\maketitle
\tableofcontents
\section{Introduction} \label{sec: Introduction}
Let $(X,\|\cdot\|)$ be a normed space. In the following
definition 
$d(x,A)=\inf\{\|x-a\|: a \in A\}$ denotes the
distance from $x$ to the set $A$.
\begin{definition}
A set $A \subseteq X$ is \textit{approximately convex} if 
$$d(tx+(1-t)y,A) \le 1$$
for all $x,y \in A$ and $t \in [0,1]$.
\end{definition} 
Recall that the Hausdorff distance between subsets $A$ and $B$
of $X$ is defined by
$$ \mathcal{H}(A,B) = \sup\{d(x,B), d(y,A): x \in A, y\in B\}.$$ Thus, $A$
is approximately convex if and only if
$$ \sup_{t \in [0,1]}\mathcal{H}(tA+(1-t)A,A) \le 1.$$
The aim of this article is to study the
relationship betwen the 
 size of an approximately convex set, as measured by its \textit{diameter}
$$\operatorname{diam}(A) = \sup\{\|x-y\|: x,y \in A\},$$
and the extent to which $A$ fails to be convex, as measured by 
the Hausdorff distance $\mathcal{H}(A,\Co(A))$ from $A$ to
its \textit{convex hull} $\Co(A)$. 

In Section~\ref{sec: appconvexsets} we extend some of the results of
\cite{Di-Ho-Ro-1} to the case of approximately convex sets. In particular,
it is shown that if $X$ is an $n$-dimensional normed space then the quantity
$$ C(X) = \sup\{\mathcal{H}(A,\Co(A)): \text{$A\subseteq X$ is approximately convex}\}$$
satisfies \begin{equation} \label{C(X)bounds}
\log_2n \le C(X) \le \lceil \log_2(n+1) \rceil, \end{equation}
where $\lceil x \rceil$ denotes the smallest integer $n \ge x$.
For the Euclidean spaces $\mathbb{R}^n$, we prove that $C(\mathbb{R}^n)=\log_2n$
for infinitely many values of $n$. Thus, the lower bound in \eqref{C(X)bounds}
is sharp. 

We also prove in Section~\ref{sec: appconvexsets} that \textit{every}
infinite-dimensional normed space contains an approximately convex set $A$
with $\mathcal{H}(A,\Co(A))=\infty$. This is used to show that the Hyers-Ulam
stability theorem fails rather spectacularly in every infinite-dimensional normed space.

In our previous paper \cite{Di-Ho-Ro-1} we studied the quantity 
$\mathcal{H}(A,\Co(A))$ 
for 
the class of \textit{approximately Jensen-convex sets}  defined as follows.
\begin{definition}
A set $A \subseteq X$ is \textit{approximately Jensen-convex} if 
$$d\left(\frac{x+y}{2},A\right) \le 1$$
for all $x,y \in A$.
\end{definition}

Suppose again that $X$ is an $n$-dimensional normed space.
In the construction of approximately convex sets $A\subseteq X$ 
presented in Section~\ref{sec: appconvexsets}, we find  that 
$\operatorname{diam}(A) \rightarrow \infty$
as $\mathcal{H}(A,\Co(A))$ approaches $C(X)$. Section~\ref{sec: diamupperbounds}
refines this construction to produce such
sets whose diameters are not too large 
in an asymptotic sense as 
$n \rightarrow \infty$. To make this precise, let us say that 
an approximately convex set $A$ is \textit{bad} if 
$\mathcal{H}(A,\Co(A)) \ge \log_2n -1$. Then our main result says that
every $n$-dimensional normed space contains bad approximately
convex sets of diameter $O(\sqrt n (\log n)^2)$. The proof uses
a result of Bourgain and Szarek \cite{Bo-Sz-1} from the local theory of Banach spaces.

In Section~\ref{sec: diameuclid} we show that the factor $\sqrt n$ in the latter result 
is sharp by demonstrating the \textit{lower} bound $\operatorname{diam}(A) \ge 0.76\sqrt{n}$
for \textit{all} bad approximately convex  sets in the Euclidean space $\mathbb{R}^n$
when $n$ is sufficiently large. We also construct nearly extremal approximately convex sets in
$\mathbb{R}^n$ of diameter $O(\sqrt{n\log n})$, which is better than our estimate 
in the general normed
space case.

Our constructions  uses the clasical \textit{entropy function}
$$E_n(t_1,\dots,t_{n+1})=\sum_{i=1}^{n+1} t_i\log_2(1/t_i)$$ 
defined on the standard $n$-simplex. In particular, we make heavy use of the fact 
that $E_n$ is an 
\textit{approximately convex function}. 
This observation seems to be new, and we include its short proof in 
Section~\ref{sec: appconvexfunctions}. 
As a corollary
we obtain the \textit{best constants} in the classical 
Hyers-Ulam stability theorem \cite{Hy-Ul-1} when $n+1$ is a power of $2$.

The last two sections concern approximately convex sets 
in infinite-dimensional spaces. The results in
 Section~\ref{sec: diamlowertype} are in principle not new: they are
essentially reformulations of known results of Larsson \cite{Lar-1} and
of Casini and Papini \cite{Ca-Pa-1} (also of Bruck \cite{Br-1}). 
It is shown that 
$X$ is $B$-convex if and only if
there exists $c >0$ such that 
$$\operatorname{diam}(A) \ge c \exp(c\mathcal{H}(A,\Co(A)))$$
for every approximately convex set $A\subseteq X$. A  similar
bound with a sharp exponent is given for spaces of type $p$.

Our deepest and perhaps most interesting result
is Theorem~\ref{thm: worst1} of Section~\ref{sec: infinitedim},
which says that the trivial inequality
$$\operatorname{diam}(A)\ge \mathcal{H}(A,\Co(A))$$ 
is actually \textit{best possible}
in general Banach spaces. More precisely, we show that for every $M>0$ there exists a 
Banach space $X$ (which is isomorphic to $\ell_1$) and an approximately convex set 
$A \subseteq X$ such that
$$\operatorname{diam}(A)\ = \mathcal{H}(A,\Co(A))=M.$$ 
The space $X$ is obtained from a rather complicated combinatorial construction 
which may conceivably have other applications in Banach space theory.
Theorem~\ref{thm: worst1} and its proof may be read independently of the rest of the paper.

Finally, a few words about notation. 
All normed spaces are assumed to be \textit{real}.
The closed unit ball $\{x \in X: \|x\|\le1\}$
of a normed space $X$ is denoted  $B_X$.
The closed ball of radius $R$ is
denoted  $B_R(X)$. The dual space of $X$ is denoted  $X^*$. A closed subspace $Y$
of $X$ has a \textit{finite-dimensional decomposition} if there exist
finite-dimensional subspaces $F_n \subseteq Y$ ($n \ge 1$) such that every $y \in Y$
admits a unique representation as a convergent series $y = \sum_{n=1}^\infty y_n$
with $y_n \in F_n$. This implies that the
finite-dimensional projections $P_n(y) = \sum_{i=1}^n y_i$ are uniformly
bounded
in the operator norm.
We write $Y=\sum_{n=1}^\infty \oplus F_n$. The
sequence spaces $\ell_p$, the finite-dimensional spaces $\ell_p^n$,
the Lebesgue spaces $L_p(0,1)$  ($1 \le p \le \infty$), and the space $C(0,1)$
of all continuous functions on $[0,1]$,
are all equipped with their classical norms. 
More specialized terminology from 
Banach space theory will be introduced as  needed.

\section{Approximately convex functions} \label{sec: appconvexfunctions}
Hyers and Ulam \cite{Hy-Ul-1} introduced the notion of an $\varepsilon$-convex function.
\begin{definition} \label{def: epsilonconvexfunc}
 Let $C$ be a convex subset of $X$ and let $\varepsilon\ge0$.
A function $f:C \rightarrow \mathbb{R}$ is $\varepsilon$-convex if
\begin{equation}
f(tx+(1-t)y) \le tf(x)+(1-t)f(y) + \varepsilon
\end{equation} for all $x,y \in C$ and $t \in [0,1]$.
\end{definition}
Note that if $f$ is $\varepsilon$-convex then  the function
$\lambda f$ is $\lambda\varepsilon$-convex for each $\lambda >0$.
Thus, $\varepsilon$ merely plays the role of a scaling factor. 
For our results it is convenient to normalize by taking $\varepsilon=1$
as follows.
\begin{definition}\label{def: appconvexfunc} Let $C$ be a convex subset of $X$.
A function $f:C \rightarrow \mathbb{R}$ is \textit{approximately convex} if
\begin{equation}
f(tx+(1-t)y) \le tf(x)+(1-t)f(y) + 1
\end{equation} for all $x,y \in C$ and $t \in [0,1]$.
\end{definition}

For $n \ge 1$, let
$\Delta_n=\{t=(t_i)_{i=1}^{n+1}:t_i \ge 0, \sum_{i=1}^{n+1}t_i=1\}$
be the standard
$n$-simplex. Let $e_i$ ($1 \le i \le n+1$) be the vertices of $\Delta_n$ and
let $\mathcal{F}_n$ be the collection of all approximately
convex functions $f: \Delta_n \rightarrow \mathbb{R}$ satisfying $f(e_i) \le 0$
for $1 \le i \le n+1$. Now define \begin{equation}
\kappa(n)= \sup_{f \in \mathcal{F}_n} \sup_{x \in \Delta_n} f(x).
\end{equation}

Cholewa \cite{Ch-1} (cf.\ \cite{Hy-Is-Ra-1}) proved the following sharp version 
of the famous Hyers-Ulam stability theorem \cite{Hy-Ul-1}.
\begin{theoremA} \cite{Ch-1} Let $U\subseteq\mathbb{R}^n$ be a convex set 
and let $\varepsilon>0$.
For every $\varepsilon$-convex function
$f:U \rightarrow \mathbb{R}$ 
there exist convex functions $g$
and $g_0$ such that \begin{equation} \label{eq: hyersulam}
f(x) \le g(x) \le f(x)+\kappa(n) \varepsilon \qquad \text{and}\qquad
|f(x)-g_0(x)| \le \frac{\kappa(n)}{2} \varepsilon. \end{equation}
Moreover, $\kappa(n)$ is the sharp constant in \eqref{eq: hyersulam}
and satisfies the upper bound $\kappa(n) \le k$ for $2^{k-1} \le n < 2^k$,
i.e.\ $\kappa(n) \le \lceil \log_2(n+1) \rceil$. 
\end{theoremA}
\begin{rem} Lazckovich \cite{La-1} observed that $\kappa(n) $ is the sharp
constant for \textit{every} convex $U$ with nonempty interior. \end{rem}
The following lemma will be used repeatedly.
\begin{lem} \label{lem: convex-k(n)} Let
$f:C \rightarrow \mathbb{R}$ be approximately convex, where $C\subset X$ is convex.
Suppose that $n \ge 1$ and that $x_1,\dots,x_{n+1} \in C$. Then
\begin{equation} \label{eq: convex-k(n)}
f\left(\sum_{i=1}^{n+1}t_ix_i\right) \le \sum_{i=1}^{n+1}t_if(x_i) + \kappa(n) \end{equation}
for all $(t_i)_{i=1}^{n+1} \in \Delta_n$. \end{lem}
\begin{proof} Define $F:\Delta_n \rightarrow \mathbb{R}$ by
$$ F(t) = f\left(\sum_{i=1}^{n+1}t_ix_i\right) - \sum_{i=1}^{n+1} t_if(x_i).$$
Then $F$ is approximately convex and $F(e_i)=0$  for $1 \le i \le n+1$.
So $F(t) \le \kappa(n)$ for all $t \in \Delta_n$, which gives 
\eqref{eq: convex-k(n)}. \end{proof} 

For
our results on approximately convex sets we require a good 
\textit{lower bound} for $\kappa(n)$: we shall
show  that $\kappa(n) \ge \log_2(n+1)$, which improves the bound
$\kappa(n) \ge (1/2)\log_2(n+1)$ given in \cite{La-1}. 
 
We require the following lemma from \cite{Di-Ho-Ro-1} concerning the function $\phi(t)$
defined by $\phi(0)=0$ and
$\phi(t)=-t\log_2 t$ ($t \in (0,1]$). For completeness we include the proof.
\begin{lem} \label{lem: phi} For all $t,x,y \in [0,1]$, we have
$$ 0 \le \phi(tx+(1-t)y) -t\phi(x)-(1-t)\phi(y)\le \phi(t)x+ \phi(1-t)y.$$
\end{lem}
\begin{proof} The left-hand inequality just  says that $\phi$ is concave (to see this note
that
$\phi''(t)= -1/(t\ln 2) <0$). To prove the right-hand inequality,
first consider the case 
$0<x \le y\le1$.  For fixed $t$ and $y$, let
$$ \psi(x) = \phi(tx+(1-t)y)-t\phi(x)-(1-t)\phi(y).$$
Then
$$\psi'(x) = \frac{t}{\ln 2}(\ln x - \ln(tx + (1-t)y)) \le 0.$$
Thus $\psi(x)$ is decreasing on $[0,y]$ and  attains its maximum at $x=0$.
But \begin{align*}
\psi(0)&= \phi((1-t)y)-(1-t)\phi(y)\\
&=-(1-t)y\log_2((1-t)y)+(1-t)y\log_2y\\
&= -(1-t)y\log_2(1-t) = \phi(1-t)y. \end{align*} Thus, if $x \le y$, then
$$
\phi(tx+(1-t)y) -t\phi(x)-(1-t)\phi(y)\le \phi(1-t)y \le \phi(t)x+ \phi(1-t)y.
$$
Similarly, if $y \le x$, then $$ 
\phi(tx+(1-t)y) -t\phi(x)-(1-t)\phi(y)\le \phi(t)x \le \phi(t)x+ \phi(1-t)y.$$
\end{proof}
The approximately convex sets which we construct
in the next section are essentially graphs of the \textit{entropy} 
functions 
$$E_n(t_1,\dots,t_{n+1})= \sum_{i=1}^{n+1} t_i \log_2(1/t_i) \qquad
 ((t_i)_{i=1}^{n+1} \in \Delta_n).$$
The following crucial observation seems to be new.
\begin{prop} $E_n$ is a continuous
concave approximately convex function on $\Delta_n$.
In particular, $E_n$ is \textit{approximately affine}, i.e.\
\begin{equation} \label{eq: E_nappaffine}
|E_n(tx+(1-t)y)-tE_n(x)-(1-t)E_n(y)| \le 1 \end{equation}
for all $x,y \in \Delta_n$ and $t \in [0,1]$. 
\end{prop}
\begin{proof} $E_n(t)=\sum_{i=1}^{n+1} \phi(t_i)$ is a sum of concave functions (by
\linebreak
Lemma~\ref{lem: phi}) and so $E_n$ is concave.
For $x=(x_i)_{i=1}^{n+1}$ and $y=(y_i)_{i=1}^{n+1}$ in $\Delta_n$
and $t \in [0,1]$, we can use Lemma~\ref{lem: phi} for the first inequality
to get
\begin{align*}
E_n(tx + (1-t)y)-&tE_n(x)-(1-t)E_n(y)
\\
&=\sum_{i=1}^{n+1}(\phi(tx_i+(1-t)y_i)-t\phi(x_i)-(1-t)\phi(y_i))\\
&\le  \sum_{i=1}^{n+1}(\phi(t)x_i+ \phi(1-t)y_i)\\
&=\phi(t)\sum_{i=1}^{n+1}x_i + \phi(1-t)\sum_{i=1}^{n+1}y_i\\
&=\phi(t)+\phi(1-t).\end{align*}
The function $\phi(t)+\phi(1-t)$ is concave and symmetric about $t=1/2$.
Thus,
$$\phi(t)+\phi(1-t) \le 2\phi(1/2) =1,$$
with equality in the last inequality only if $t=1/2$. \end{proof}
\begin{rem} The fact that $E_n$
has the weaker property of 
being  approximately \textit{Jensen-convex}
 (which corresponds to setting $t=1/2$ in Definition~\ref{def: appconvexfunc})
is well-known and has been observed by various authors, e.g.\ \cite{La-1}.
\end{rem}
Note that the following theorem gives the sharp constant in the Hyers-Ulam stability theorem
when $n+1$ is a power of $2$.
\begin{thm} The constants $\kappa(n)$ satisfy the bounds
\begin{equation}
\log_2(n+1) \le \kappa(n) \le \lceil \log_2(n+1) \rceil. \end{equation}
In particular, $\kappa(n) = \log_2(n+1)$ when $n+1$ is a power of $2$.
\end{thm}
\begin{proof} The upper bound is due to Cholewa \cite{Ch-1}. For the lower bound,
since $E_n \in \mathcal{F}_n$, we have
$$ \kappa(n) \ge \max_{t \in \Delta_n} E_n(t) = E_n(1/(n+1),\dots,1/(n+1))
= \log_2(n+1).$$
\end{proof}
\begin{rem} Obviously, $\kappa(1)=1$. Green  \cite{Gr-1} showed that $\kappa(2) = 5/3$.  
In a later paper we shall show that, for $n \ge 1$,
$$\kappa(n)= [\log_2(n+1)]+2-\frac{2^{1+[\log_2(n+1)]}}{n+1},$$
where $[x]$ is the greatest integer function. 
The proof is too long to
be included  here. The corresponding constants for bounded Jensen-convex functions were
computed in \cite{Di-Ho-Ro-1}.\end{rem}

\section{Approximately convex sets} \label{sec: appconvexsets}
\begin{thm}\label{thm: C(X)} let $X$ be an $n$-dimensional normed space. There is a least
positive constant $C(X)$ such that 
\begin{equation} \label{eq: C(X)def}
\mathcal{H}(A,\Co(A)) \le C(X) \sup_{t \in [0,1]}\mathcal{H}(tA+(1-t)A,A)
\end{equation} for every nonempty $A \subseteq X$.
Moreover, $C(X)$ satisfies \begin{equation}
\log_2n  \le C(X) \le \kappa(n). \end{equation}
In particular, $\log_2n \le C(X) \le \lceil \log_2(n+1)\rceil\le \log_2n + 1$.
\end{thm}
\begin{proof} We may assume that the right-hand side of \eqref{eq: C(X)def} is
finite, otherwise there is nothing to prove.
Observe that the effect of replacing $A$ by $\lambda A$ is
to multiply
 both sides of \eqref{eq: C(X)def} by $|\lambda|$. So, by choosing $\lambda$
appropriately, we may
assume  that 
$$\sup_{t \in [0,1]}\mathcal{H}(tA+(1-t)A,A)=1.$$
The right-hand
estimate for $C(X)$
 is due to Casini and Papini \cite{Ca-Pa-1}. For completeness we
recall the proof. Let $f(x) = d(x,A)$ ($x \in X$). 
First note that $f$ is $1$-Lipschitz and non-negative. To see that
$f$ is
approximately convex, 
note that for $x,y \in X$, $a,b \in A$,  and $t \in [0,1]$, we have
\begin{align*} f(tx+(1-t)y)&=d(tx+(1-t)y,A) \\
&\le \|(tx+(1-t)y)- (ta+(1-t)b)\| \\&+ d(ta+(1-t)b,A)\\
&\le t\|x-a\| +(1-t)\|(y-b)\|+1. \end{align*}
Taking the infimum of this expression over all choices of $a$ and $b$ yields
$$f(tx+(1-t)y)\le tf(x) +(1-t)f(y) +1.$$ Now suppose that 
$x \in \Co(A)$. By Carath\'eodory's Theorem (see e.g.\ \cite[Thm.\ 17.1]{Ro-1}), 
$x=\sum_{i=1}^{n+1} t_ia_i$, a convex combination of $n+1$ 
elements $a_i \in A$.  Then Lemma~\ref{lem: convex-k(n)} yields
$$f(x) \le \sum t_if(a_i) +\kappa(n)=\kappa(n),$$
since $f(a)=0$ for all $a \in A$. 
The left-hand inequality uses the entropy functions $E_n$.   
Let
$(e_i)_{i=0}^{n-1}$ be an Auerbach basis for $X$ (see e.g.\ \cite[p. 16]{Li-Tz-1}). Recall
that this means that \begin{equation} \label{eq: Auerbach}
\max |a_i| \le \left\|\sum_{i=0}^{n-1} a_ie_i\right\| \le \sum_{i=0}^{n-1} |a_i|.
\end{equation} 
for all scalars $a_0,\dots,a_{n-1}$. Set $e_n=0$ so that $\Co\{e_i: 1 \le i \le n\}$
is an $(n-1)$-simplex.
For each $M >0$, we define a set $A_M$ thus:
\begin{equation*}
A_M = \left\{M \sum_{i=1}^{n-1} t_ie_i + E_{n-1}(t_1,\dots,t_n)e_0: 
(t_i)_{i=1}^{n} \in \Delta_{n-1}\right\}
\end{equation*}
First let us verify that $A_M$ is approximately convex. Suppose
that $0 \le t \le1$ and that
 $a=M \sum_{i=1}^{n-1} x_ie_i + E_{n-1}(x)e_0$
and $b=M \sum_{i=1}^{n-1} y_ie_i + E_{n-1}(y)e_0$
belong to $A_M$, where $x=(x_i)_{i=1}^{n}$ and $y=(y_i)_{i=1}^n$ belong
to $\Delta_{n-1}$. Then 
$c = M \sum_{i=1}^{n-1} z_ie_i + E_{n-1}(z)e_0$ 
also belongs to $A_M$, where $z=tx+(1-t)y$.
Since $e_0$ is a unit vector and $E_{n-1}$ 
is \textit{approximately affine} \eqref{eq: E_nappaffine}, we have \begin{align*}
\|ta+(1-t)&b-c\|  \\
&=|tE_{n-1}(x)+(1-t)E_{n-1}(y)
-E_{n-1}(tx+(1-t)y)| \\
&\le 1,
\end{align*} and so $A_M$ is approximately convex.
Note that $x_0= (M/n) \sum_{i=1}^{n-1}e_i \in \Co(A_M).$
We shall show that $d(x_0,A_M) \rightarrow \log_2n$ as $M \rightarrow \infty$.
To see this, fix $\varepsilon >0$. By continuity of $E_{n-1}$ there exists
$\alpha>0$ such that if $\max_{1 \le i \le n-1} |t_i - 1/n| \le \alpha$
then $E_{n-1}(t_1,\dots,t_{n}) \ge \log_2n-\varepsilon$, whence
by \eqref{eq: Auerbach}
\begin{align*}
\left\|x_0-\left(M\sum_{i=1}^{n-1}t_ie_i+E_{n-1}(t_1,\dots,t_{n})e_0\right)\right\|
&\ge E_{n-1}(t_1,\dots,t_{n})\\ &\ge \log_2n-\varepsilon. \end{align*}
Now suppose, on the other hand, that $\max_{1\le i \le n-1} |t_i - 1/n| \ge \alpha$.
By  \eqref{eq: Auerbach}
\begin{align*}
\left\|x_0-\left(M\sum_{i=1}^{n-1}t_ie_i+E_{n-1}(t_1,\dots,t_{n})e_0\right)\right\|
&\ge M\max_{1\le i \le n-1} |t_i - 1/n|\\ &\ge M\alpha \rightarrow\infty\end{align*}
as $M \rightarrow \infty$.
Thus, for all sufficiently large $M$, we have
$d(x_0, A_M) \ge \log_2n -\varepsilon$. Since $\varepsilon>0$ is arbitrary, this 
gives the lower bound $C(X) \ge \log_2n$.
\end{proof}
For large $n$ the lower bound $C(X) \ge \log_2n$ is 
actually attained for
certain Euclidean spaces (e.g.\ for $X=\mathbb{R}^{16}$).
\begin{thm} \label{thm: power2}
Suppose that $n=2^k$, where $k\ge4$. Then $C(\mathbb{R}^n)=\log_2n$.
\end{thm}
\begin{proof} For $n=2^k$, we have $\kappa(n-1)= \log_2n$. The argument
used to prove Theorem~3.7 of 
\cite{Di-Ho-Ro-1} (too lengthy to recall here)
 shows that the result will follow provided $n=2^k$ is large enough to
ensure that
$$ \kappa(n-1) \ge \frac{\sqrt{2n}(\sqrt{2n}+\sqrt{n-1})}{n+1}.$$
This holds for $k\ge 4$. \end{proof}
\begin{rem} The calculation of $C(\mathbb{R}^n)$ for small $n$ seems
problematic. Clearly $C(\mathbb{R})=1$, and
examples show that $C(\mathbb{R}^2) > 1.37$.
In \cite{Di-Ho-Ro-1} the 
corresponding constants for approximately Jensen-convex sets in $\mathbb{R}^n$
were computed
in all dimensions. \end{rem}

Before turning to infinite-dimensional spaces, let us 
make the following definition (the analogue
of Definition~\ref{def: epsilonconvexfunc}). 
\begin{definition} \label{def: epsilonconvexset}
Let $\varepsilon>0$. A set $A \subseteq X$ is $\varepsilon$-convex
if $$d(ta+(1-t)b,A)\le \varepsilon$$ for all $a,b \in A$.
\end{definition}

\begin{thm} \label{thm: infdimset}
 Let $X$ be an infinite-dimensional normed space. There exists
an approximately convex set $A \subseteq X$ such that $\mathcal{H}(A,\Co(A))=\infty$.\end{thm}
\begin{proof} We shall use the following consequence of Theorem~\ref{thm: C(X)}.
Let $\varepsilon>0$ and $M>0$. Then every normed space of sufficiently
large dimension contains a compact $\varepsilon$-convex set $A$ such
that $\mathcal{H}(A,\Co(A)) > M$. Using this fact repeatedly, a
routine argument (cf.\ \cite[p. 4]{Li-Tz-1}) shows that $X$ contains a
subspace $Y$ with a \textit{finite-dimensional decomposition}
$\sum_{n=1}^\infty \oplus F_n$ and sets $A_n\subseteq F_n$ ($n \ge 1$)
such that $A_n$ is a $2^{-n}$-convex set containing zero and
$\mathcal{H}(A_n,\Co(A_n)) > n$. Let $A$ be the collection of all
vectors of the form $\sum_n x_n$, where $x_n \in A_n$ and only
finitely many of the $x_n$'s are nonzero.

 First let us verify that $A$ is approximately convex.
Suppose that $x=\sum_n x_n$ and $y = \sum_n y_n$ are in $A$ and that $0 \le t \le 1$.
Since $A_n$ is $2^{-n}$-convex and compact, there exists $z_n \in A_n$ with
$\|z_n-(tx_n+(1-t)y_n)\| \le 2^{-n}$. Moreover, we may choose the $z_n$'s so that only 
finitely many are nonzero, ensuring that
$z = \sum_n z_n$ belongs to $A$. By the triangle inequality
\begin{equation*}
\|z -(tx+(1-t)y)\| \le \sum_n \|z_n-(tx_n+(1-t)y_n)\| \le \sum_n 2^{-n} = 1.
\end{equation*}

Let us verify that $\mathcal{H}(A,\Co(A))=\infty$.
Since $\sum_{n=1}^\infty \oplus F_n$ is a 
finite-dimensional decomposition,
the natural
projection maps
from  \linebreak $\sum_{n=1}^\infty \oplus F_n$ onto $F_n$ are uniformly bounded 
in operator norm
by $K$, say.
Since $\mathcal{H}(A_n,\Co(A_n)) > n$,
there exists $w_n \in \Co(A_n)$ such that $d(w_n,A_n) \ge n$, and 
since $A_n \subseteq F_n$, we have $$d(w_n,A) \ge (1/K)d(w_n,A_n) \ge n/K.$$
Thus, $\mathcal{H}(A,\Co(A))=\infty$.
\end{proof}
As an application of the last result we show that the Hyers-Ulam stability
theorem (Theorem~A above) 
 fails rather dramatically in \textit{every} infinite-dimensional normed space
(cf.\ \cite{Ca-Pa-1}).
\begin{cor} \label{cor: HUinfdim}
 Let $X$ be an infinite-dimensional normed space. There exists a 1-Lipschitz
approximately convex function $f:X \rightarrow \mathbb{R}$ with the following property.
For all $M>0$ there exists $R>0$ such that for every convex function $g: B_R(X) \rightarrow
\mathbb{R}$, we have
\begin{equation*}
\sup_{x \in B_R(X)} |f(x) - g(x)| > M.  \end{equation*}
In particular, $\sup\{|f(x)-g(x)|: x \in X\}=\infty$ for 
every convex function $g:X \rightarrow \mathbb{R}$.
\end{cor}
\begin{proof} Using the notation of Theorem~\ref{thm: infdimset},
 we prove that 
$f(x) = d(x,A)$ has the required property.
It was shown in Theorem~\ref{thm: C(X)} that
$f$ is approximately convex and $1$-Lipschitz.
Choose $R$ so that
$\Co(A_n) \subseteq B_R(X)$. Suppose that $g: B_R(X)\rightarrow \mathbb{R}$ is a convex
function satisfying
$|g(x) - f(x)| \le M$. Since $f(x)=0$ for all $x \in A_n$, it follows that
$g(x) \le M$ for all $x \in A_n$, and hence $g(x) \le M$ for all $x \in \Co(A_n)$.
But $f(w_n) > n$, and so $M > n/2$.  \end{proof}

Recall that a normed space $X$ is \textit{B-convex}
 if $X$ does not `contain $\ell_1^n$'s
uniformly', i.e., if there exist $n \ge 2$ and $\alpha>0$ such that 
$$ \min_{\pm}\left\| \sum_{i=1}^n \pm x_i\right\| \le n-\alpha$$
for all $x_i \in B(X)$ ($1 \le i \le n$).

For general normed spaces, Corollary~\ref{cor: HUinfdim} is close to optimal
in view of the following positive result on the approximation
of Lipschitz $\varepsilon$-convex functions on \textit{bounded} sets
from \cite{Ca-Pa-2}.
(Here (a)$\Rightarrow$(c) is \cite[Thm.\ 1]{Ca-Pa-2} and
(b)$\Rightarrow$(a) is
implicit in  \cite[Props.~1,2]{Ca-Pa-2}. The other
implication (c)$\Rightarrow$(b) is  trivial.)

\begin{theoremB} \cite{Ca-Pa-2} Let $X$ be a normed space. The following are equivalent:
\begin{itemize}
\item[(a)] $X$ is B-convex;
\item[(b)] there exist $k < 1/2$ and $\alpha>0$ such that for every
$\varepsilon<\alpha$ and for every
$\varepsilon$-convex $1$-Lipschitz function $f:B(X)\rightarrow\mathbb{R}$ there exists a
convex function $g: B(X) \rightarrow \mathbb{R}$ such that
\begin{equation*} |g(x)-f(x)| \le k \qquad(x \in B(X));
\end{equation*}
\item[(c)] there exist $c>0$ and $\alpha>0$
 such that for every $\varepsilon < \alpha$ and for every
$\varepsilon$-convex $1$-Lipschitz function $f:B(X)\rightarrow\mathbb{R}$ there exists a
convex function $g: B(X) \rightarrow \mathbb{R}$ such that
\begin{equation*} |g(x)-f(x)| \le c\varepsilon\log_2(1/\varepsilon)\qquad(x \in B(X)).
\end{equation*} \end{itemize}\end{theoremB}

\begin{rem} \label{rem: approxlip} Condition (b) of this result is very
pertinent to
Section~\ref{sec: infinitedim} below, where  we prove (Corollary~\ref{cor: approxlip}) that 
for $X=C(0,1)$ there is no constant $k<1$  such
that (b) holds. This is clearly an optimal result
since every $1$-Lipschitz function $f$ on $B(X)$ 
satisfies $|f(x)- c| \le 1$, where $c= (\inf f+ \sup f)/2$, i.e.\ (b) holds for $k=1$.
\end{rem}
\section{Diameter of approximately convex sets} \label{sec: diamupperbounds}
Our next goal is to prove
that {\it every} $n$-dimensional normed space contains  a ``bad''
approximately convex
set  (that is, $\mathcal{H}(A,\Co(A)) 
\ge \log_2(n+1) -\varepsilon$)
of diameter $O(\sqrt n (\log n)^2)$.
In the next section we shall prove that for Euclidean spaces this estimate 
for the diameter is fairly sharp.

For two isomorphic Banach spaces $X$ and $Y$  
recall
that their Banach-Mazur distance $d(X,Y)$ is defined thus:
$$ d(X,Y)= \inf\{\|T\|\|T^{-1}\|:\text{$T:X \rightarrow Y$ is an isomorphism}\}.$$
\begin{thm} \label{th: upbound1}
Let $\varepsilon \in (0,3)$. For all sufficiently large $n$ and
all
normed spaces $X$ of dimension $n$
there exists an approximately convex set
$A \subseteq X$ such that
\begin{equation}\label{eq: firstcon} 
\mathcal{H}(A, \Co(A)) \ge \log_2n -\varepsilon
\end{equation} and
\begin{equation} \label{eq: secondcon}
\operatorname{diam}(A)\le \frac{25}{\varepsilon}(\log_2n)^2 d(X,\ell_1^n).
\end{equation}
\end{thm}
\begin{proof} In order to simplify notation we shall prove the result for all normed
spaces $X$ of dimension $n+1$ (with $n+1$ replacing $n$ in \eqref{eq: firstcon}
and \eqref{eq: secondcon}). Note that $X$ contains a subspace $Z$ of codimension one
such that $d(Z,\ell_1^n) \le d(X,\ell_1^{n+1})$ (since $\ell_1^{n+1}$ contains 
subspaces isometric to
$\ell_1^n$). Let $F$ be a linear functional in $X^*$ of unit norm such that
$Z=\operatorname{ker}(F)$. Let
$e_0$ be an unit vector in $X$ which is normed by $F$, i.e.,
such that $F(e_0)=\|e_0\|=1$.
Note that by the triangle inequality
\begin{equation} \label{eq: triangle}
 \|z+ \lambda e_0\| \ge 
\max(\|z\|-|\lambda|, |\lambda|) \ge\max\left(\frac{\|z\|}{2},|\lambda|\right)
\end{equation}
for all $z \in Z$ ($=\operatorname{ker}(F)$) and $\lambda \in \mathbb{R}$. 
Since $d(Z,\ell_1^n) \le d(X,\ell_1^{n+1})$,
$Z$ has a basis $(e_k)_{k=1}^n$ satisfying
\begin{equation} \label{eq: Banachmazur}
\sum_{k=1}^{n}|a_k| \le \left\|\sum_{k=1}^{n} a_k e_k\right\|
\le d(X,\ell_1^{n+1})\sum_{k=1}^{n}|a_k| ,\end{equation}
for all choices of scalars $(a_k)_{k=1}^{n}$.
For each  $M>0$, define $A_M \subseteq X$ by
$$ A_M = \left\{M\left(\sum_{k=1}^n t_k e_k\right) + 
E_{n-1}(t_1,\dots,t_n)e_0:
(t_1,\dots,t_n)\in \Delta_{n-1}\right\}.$$
It was proved in Theorem~\ref{thm: C(X)} that
$A_M$ is approximately convex for all choices of $M$. 
Observe also that 
$$x_0=\frac{M}{n}\sum_{k=1}^n e_k
\in \Co(A_M).$$ In order to verify \eqref{eq: firstcon}, it suffices to show
that $$d(x_0, A_M) \ge \log_2(n+1)-\varepsilon$$ for a suitable choice of $M$.
To that end, fix $\alpha\in(0,1)$ and fix $y=M(\sum_{k=1}^n t_ke_k) +  
\left(\sum_{k=1}^n t_k \log_2(1/t_k)\right)e_0
\in A_M$. Let $$B_1=\{k: t_k \ge (1+\alpha)/n\},\qquad
B_2= \{k:t_k<(1+\alpha)/n\},$$ and 
set $\mu(B_i)=\sum_{k \in B_i} t_k$ ($i=1,2$). Then,
by \eqref{eq: triangle} for the first inequality and 
the left-hand side of \eqref{eq: Banachmazur} for the second, we have
\begin{align*}
 \|y-x_0\|&=\left\|\sum_1^n M
 (t_k-(1/n))e_k + \left(\sum_{k=1}^n t_k \log_2(1/t_k)\right)e_0\right\|\\
&\ge \max\left(\frac{1}{2}\left\|\sum_1^n M
 (t_k-(1/n))e_k\right\|,\sum_{k=1}^n t_k \log_2(1/t_k)\right)\\ 
 &\ge\max\left(\frac{M}{2}\sum_1^n 
 |t_k-(1/n)|,\sum_{k=1}^n t_k \log_2(1/t_k)\right) \\
&\ge \max\left(\sum_{k \in B_1}
\frac{M}{2}|t_k-(1/n)|,\sum_{k \in B_2} t_k|\log_2(t_k)|\right)\\
&\ge \max\left(  \frac{M}{2}\left(\frac{\alpha}{2}\sum_{k\in B_1}t_k\right),
\left(\sum_{k\in B_2}t_k\right)
(\log_2n-\log_2(1+\alpha))\right)\\
\intertext{(since $t_k - (1/n) \ge (\alpha/(1+\alpha))t_k \ge (\alpha/2)t_k$
for $k \in B_1$)}
&\ge \max\left(\frac{M\alpha}{4}\mu(B_1), (\log_2n-(3\alpha/2))\mu(B_2)\right), 
\end{align*} where at the last step we use the fact that $\log_2(1+\alpha)
\le 3\alpha/2$ for $\alpha \in [0,1]$.
Now set $M=4(\log_2n)^2/\alpha$. There are two cases to consider. First, if 
$\mu(B_2) \ge 1- \alpha/\log_2n$, then \begin{multline*}
\|y-x_0\| \ge \left(\log_2n-\frac{3\alpha}{2}\right)\mu(B_2)\\ 
\ge 
\left(\log_2n-\frac{3\alpha}{2}\right)(1-\frac{\alpha}{\log_2n})
 \ge \log_2n-\frac{5\alpha}{2}.\end{multline*}
Secondly, if $\mu(B_1) \ge \alpha/\log_2n$, then
$$ \|y-x_0\| \ge \frac{M}{4}\frac{\alpha}{\log_2n} = \log_2n.$$
Hence $\|y-x_0\| \ge \log_2n-(5\alpha/2)$. Setting $\alpha=\varepsilon/3$
we see that \eqref{eq: firstcon}
is satisfied (with $n$ replaced by $n+1$) by $A=A_M$ whenever
$n$ is large enough to ensure that $\log_2(n+1) - \log_2(n) \le \alpha/2$.
Finally, the right-hand side of \eqref{eq: Banachmazur} yields
\begin{align*}
 \operatorname{diam}(A) &\le 2d(X,\ell_1^{n+1})M+ \log_2n\\
&\le d(X,\ell_1^{n+1})(4(\log_2n)^2/\alpha) + \log_2n\\
&\le 25(\log_2n)^2d(X,\ell_1^{n+1})/\varepsilon, \end{align*}
for all sufficiently large $n$, and so $A$ satisfies condition \eqref{eq: secondcon}. 
\end{proof}
Since $d(\ell_p^n,\ell_1^n)=n^{(p-1)/p}$ for $1 \le p \le 2$, we get the following corollary.
\begin{cor} \label{cor: l_pcase} Let $1 < p \le 2$ and let $\varepsilon \in(0,3)$. 
For all sufficiently
large $n$ there exists an approximately convex set $A \subseteq \ell_p^n$
such that
\begin{equation*}
\mathcal{H}(A, \Co(A)) \ge \log_2n -\varepsilon \end{equation*}
and \begin{equation*}
\operatorname{diam}(A)\le \frac{25}{\varepsilon}n^{(p-1)/p}(\log_2n)^2.
\end{equation*} \end{cor}
\begin{rem} For $p=2$, a much stronger result will be proved in the next section.
\end{rem}

For $p=1$, we may reduce the exponent of $\log_2n$.
\begin{prop} \label{prop: casep=1} Let $\varepsilon \in (0,2)$. For all sufficiently large $n$
there exists an approximately convex set $A \subset \ell_1^n$ such that
\begin{equation}\label{eq: firstcon1}
\mathcal{H}(A, \Co(A)) \ge \log_2n -\varepsilon \end{equation}
and \begin{equation} \label{eq: secondcon1}
\operatorname{diam}(A)\le \left(\frac{8}{\varepsilon}+1\right)\log_2n.
\end{equation} \end{prop}
\begin{proof} Setting $X=\ell_1^{n+1}$, we 
follow Theorem~\ref{th: upbound1} taking advantage of some
simplifications in the proof which we now indicate. First,
we may choose $(e_0,\dots,e_n)$ to be the standard unit vector basis of
$\ell_1^{n+1}$, so that \eqref{eq: triangle} becomes
simply  $\|z+\lambda e_0\|= \|z\| + |\lambda|$,
for all $z= \sum_{i=1}^n a_i e_i \in Z$. The estimate for $\|y-x_0\|$ then becomes
$$\|y-x_0\|\ge \left(\frac{M\alpha}{2}\mu(B_1)+(\log_2n-(3\alpha/2))\mu(B_2)\right).$$ 
Setting $M= 2(\log_2n)/\alpha$,
we obtain
$$\|y-x_0\| \ge \left(\log_2n-\frac{3\alpha}{2}\right)(\mu(B_1)+\mu(B_2))
=\log_2n-\frac{3\alpha}{2}.$$
Setting $\alpha=\varepsilon/2$
we see that \eqref{eq: firstcon1}
is satisfied (with $n$ replaced by $n+1$) by $A=A_M$ whenever
$n$ is large enough to ensure that $\log_2(n+1) - \log_2(n) \le \alpha/2$. Finally,
$$\operatorname{diam}(A) \le 2M+\log_2n \le \left(\frac{4}{\alpha}+1\right)\log_2n,$$
which yields \eqref{eq: secondcon1}.
\end{proof}

\begin{rem}\label{rem: l1case}
 In particular, $\ell_1^n$ contains ``bad'' approximately convex sets
of ``small'' diameter $O(\log n)$. Indeed, the trivial
 lower bound 
$$\operatorname{diam}(A) \ge \mathcal{H}(A, \Co(A))
 \ge \log_2n -\varepsilon$$
shows that the diameter must grow at least logarithmically with $n$.
\end{rem}
Finally, we come to the main result of this section. 
\begin{thm} \label{th: upbound2} 
Let $\varepsilon \in (0,6)$.
For all sufficiently large $n$ and
all
normed spaces $X$ of dimension $n$
there exists an approximately convex set
$A \subseteq X$ such that
\begin{equation} \label{eq: firstcond}
\mathcal{H}(A, \Co(A)) \ge \log_2n - \varepsilon \end{equation}
and
\begin{equation} \label{eq: secondcond}\operatorname{diam}(A)\le
 \frac{K(\log_2n)^2 \sqrt n}{\varepsilon^3},\end{equation}
where $K$ is an absolute constant.
 \end{thm}
\begin{proof} Fix $\theta \in (0,1)$. Bourgain and
Szarek \cite{Bo-Sz-1} (cf.\ also \cite{Sz-1}) proved that every $n$-dimensional normed space $X$ 
contains a subspace $Y$, with $\operatorname{dim}Y=k > [\theta n]$,
satisfying \begin{equation} \label{eq: d(Y,l_1^k)}
d(Y,\ell_1^k) \le C(1-\theta)^{-2}\sqrt n, \end{equation}
where $C$ is a  constant.
Set $\theta = 1-\varepsilon/6$.  Then, for $\varepsilon<1$, \begin{equation}\label{eq: log_2k}
 \log_2k \ge \log_2n - \log_2(1/\theta) \ge \log_2n - \varepsilon/2. \end{equation}
Applying Theorem~\ref{th: upbound1}
 to $Y$ and to $\varepsilon/2$
 yields an approximately convex set
 $A \subseteq Y$
satisfying \eqref{eq: firstcond} (from \eqref{eq: log_2k} and \eqref{eq: firstcon})  
and \eqref{eq: secondcond} (from \eqref{eq: d(Y,l_1^k)} and \eqref{eq: secondcon}).
\end{proof}
\section{Bounds  in Euclidean spaces} \label{sec: diameuclid}
In this section we prove
that the ``bad'' approximately convex sets constructed in Theorem~\ref{th: upbound2}
necessarily have  diameter larger than \linebreak $0.76\sqrt n$ 
in $n$-dimensional Euclidean spaces
when $n$ is large. The proof uses only elementary  geometry. Along the way
we prove a result about Hilbert space (Theorem~\ref{th: bestconv})
 which may be of independent interest
because of its sharp constants. We also improve the upper bound
of Corollary~\ref{cor: l_pcase} by constructing a nearly extremal
approximately convex set in $\mathbb{R}^n$ of diameter $O(\sqrt{n \log n})$.

Recall that a simplex $\Sigma \subseteq \mathbb{R}^n$ is \textit{regular}
if its edges all have the same Euclidean length.
\begin{lem} \label{prop: regularbest} Let $\Sigma$ be an $n$-simplex which 
contains the origin 
in its interior and
whose vertices lie on the Euclidean unit sphere $S^{n-1}$.
For each $0 \le k \le n-1$ there exists a $k$-face
$F_k$ of $\Sigma$ such that
$$ d(0,F_k) \le \alpha_{n,k}=\sqrt{\frac{n-k}{n(k+1)}},$$
with equality if $\Sigma$ is a regular simplex. 
\end{lem}
\begin{proof} First we prove the result for $k=n-1$. Let $V$ be one of the vertices for which 
the corresponding barycentric coordinate of the origin is at most $1/(n+1)$.
Let the line segment through the origin joining $V$ to the  opposite $(n-1)$-face $F$ 
intersect $F$ in a point $P$, say.  Then the 
origin divides the line joining $V$
to $P$ into two segments bearing a ratio  of not less than $n$ to $1$. Since $V$ lies on
the unit sphere, it follows that
 $d(0,P)\le 1/n$. Thus, 
$d(0,F) \le 1/n = \alpha_{n,n-1}$, which completes the proof for the case
$k=n-1$.

 The proof for $0<k<n-1$ is by induction on $n$.
Suppose that the result holds for  $n-1$ and for $0<k<n-1$. Let $F_{n-1}$ be an
$(n-1)$-face of $\Sigma$ nearest to the origin
and let
 $Q$ be the point in $F_{n-1}$ nearest to the origin. 
Then $0 \le d=d(0,Q)= d(0,F_{n-1}) \le 1/n$.
The largest
Euclidean ball inscribed in $\Sigma$ with center the origin touches $F_{n-1}$ at $Q$. Hence
$Q$ is in the interior of the $(n-1)$-simplex $F_{n-1}$ whose vertices 
lie on the $(n-2)$-sphere with center $Q$
and radius $\sqrt{1-d^2}$. Fix $0<k<n-1$. By the inductive hypothesis applied to
$Q$ and $F_{n-1}$ there exists a $k$-face $F_k$ of  $F_{n-1}$ such that
$$d(Q,F_k) \le \alpha_{n-1,k} \sqrt{1-d^2}.$$ 
So 
$$d(0,F_k)^2 =d(0,Q)^2 + d(Q,F_k)^2 \le d^2 + (1-d^2)\alpha_{n-1,k}^2,$$
where $0 \le d \le 1/n$. The right-hand side is greatest when $d=1/n$, which
gives 
$$d(0,F_k)^2 \le \frac{1}{n^2} 
 +\left(1-\frac{1}{n^2}\right)\alpha_{n-1,k}^2 = \alpha_{n,k}^2.$$
\end{proof}
The following theorem is perhaps of independent interest because of the sharp
constants.
\begin{thm} \label{th: bestconv}
Let $(x_i)_{i=1}^{n+1}$  be elements from the unit ball 
of a Hilbert space $H$ and suppose that
 $0 \in \Co(\{x_i: 1 \le i \le n+1\})$. For each $1 \le j \le n$,
there exists $J \subseteq \{i: 1 \le i \le n+1\}$ such that $|J|=j$ and
$$ d(0, \Co(\{x_i: i \in J\})) \le \sqrt{\frac{n+1-j}{nj}}.$$ 
\end{thm}     
 \begin{proof} By slightly perturbing the elements, if necessary,
 we may assume that the set $\{x_i: 1\le i \le n+1\}$ is affinely independent
and that the origin lies in the interior of
the simplex $\Co(\{x_i: 1 \le i \le n+1\})$. 
Let $y_i= x_i/\|x_i\|$. Clearly,
$$ d(0, \Co(\{x_i: i \in A\}))
\le d(0, \Co(\{y_i: i \in A\}))$$ 
for all $A \subseteq\{i:1\le i \le n+1\}$. Now
Lemma~\ref{prop: regularbest} applied
 to the simplex $\Sigma$ with vertices $\{y_i:1\le i\le n+1\}$ yields the desired result.
  \end{proof}

\begin{thm} \label{th: lowbound3} Suppose that $A\subseteq \mathbb{R}^{\text{$n$}}$ 
is approximately convex and
satisfies $\mathcal{H}(A,\Co(A))\ge  \log_2 n -1$.
Then, for any integer $j$ with $1\le j\le n$ we have 
\begin{equation} \label{eq: diamA}
\operatorname{diam}(A) 
\ge \left( \frac{(\log_2n -1-
\lceil\log_2j\rceil)\sqrt{j}}{\sqrt{n-j+1}}\ \right)  \sqrt n
\end{equation} 
In particular, $A$ satisfies the (nontrivial) lower bounds $\diam(A)
\ge 0.7525 \sqrt n$ for all $n\ge20$, and $\diam(A) \ge 0.768 \sqrt n$
for all sufficiently large $n$.
\end{thm}
\begin{proof} 
Assuming (as we may) that $A$ is compact,
there exists $x_0\in \Co(A)$ 
with $d(x_0,A)\ge \log_2n-1$. 
By translating $A$, we may assume that
$x_0=0$. Thus, $0 \in \Co(A)$ and $d(0,A)\ge \log_2n-1$. 
The fact that $\operatorname{diam}(A)=D$ now implies
that $\|x\| \le D$ for all $x \in A$.
By Carath\'eodory's Theorem, there exist $(x_i)_{i=1}^{n+1}$
in $A$ such that $0 \in \Co(\{x_i: 1 \le i \le n+1\})$.
Let $1\le j \le n$, then by Theorem~\ref{th: bestconv} there exists
$J \subseteq \{i:1\le i \le n+1\}$ such that $|J|=j$ and
$$
d(0, \Co(\{x_i: i \in J\}))
 \le \left(\sqrt{\frac{n-j+1}{nj}}\ \right)D=
\left(\sqrt{\frac{n-j+1}{j}}\ \right)\frac{D}{\sqrt n}.
$$
Let $y_0$ be the point in $\Co(\{x_i: i \in J\})$ nearest the origin.
Because $A$ is approximately convex, the function
$d(x,A)$ is an approximately convex function which vanishes at each $x_i$.
So, by Lemma~\ref{lem: convex-k(n)}, $d(y_0,A) \le \kappa(j-1)\le
\lceil\log_2j\rceil$ for $2\le j\le n$ and if $j=1$ then $y_0\in A$
and $d(y_0,A)=0=\lceil\log_21\rceil$.  Thus $d(y_0,A)\le
\lceil\log_2j\rceil$ for $1\le j \le n$.  Therefore
$$
\log_2n -1 \le d(0,A)
\le \|y_0\| + d(y_0,A) 
\le\left(\sqrt{\frac{n-j+1}{j}}\ \right)\frac{D}{\sqrt n}+ \lceil\log_2j\rceil
$$
which yields 
$$
D \ge  \frac{\big(\log_2n -1-
\lceil\log_2j\rceil\big)\sqrt{j}}{\sqrt{n-j+1}}\,   \sqrt n = f(j,n)\sqrt{n}
$$
where this defines $f(j,n)$.  If $k$ is a non-negative integer with $2^k\le n$ then
$\lceil\log_22^k\rceil=k=\log_2(2^k)$.  Therefore
\begin{align*}
f(2^k,n)&=\frac{\big(\log_2n-1-\log_2{2^k}\big)\sqrt{2^k}}{\sqrt{n-2^k+1}}
=\frac{\log_2(n/2^k)-1}{\sqrt{(n/2^k)-1+2^{-k}}}\\
&= F(n/2^k) + r(k,n)
\end{align*}
where $F(\alpha)= (\log_2(\alpha)-1)/({\sqrt{\alpha-1}})$
and $r(k,n)\to 0$ as $n,k\to \infty$.  For each $n$ and $\alpha>0$
there is an integer $k$ so that $\alpha\le n/2^k \le 2\alpha$, and for
$\alpha_0=9.109883742$ 
$$
\alpha_0\le \alpha\le 2\alpha_0 \quad\text{implies}\quad 
	F(\alpha)\ge 0.76811996\ .
$$
Therefore if $n$ is sufficiently large and $k$ is chosen so
that $\alpha_0\le n/2^k \le 2\alpha_0$ then
$$
\max_{1\le j \le n}f(j,n)\ge f(2^k,n) \ge 0.768
$$
and thus $D\ge 0.768\sqrt{n}$.

For any $n$ and  $k\ge1$ 
$$
f(2^k,n)\ge \frac{\log_2(n/2^k)-1}{\sqrt{(n/2^k)-1+1/2^1}} = G(n/2^k)
$$
where $G(\beta) =(\log_2(\beta)-1)/\sqrt{\beta-1/2}$.  If
$\beta_0=9.919205826$, then $\beta_0\le \beta \le 2\beta_0$ implies
$G(\beta)\ge 0.7525$.  Now assume that $n\ge 20$, and that $\beta_0\le
n/2^k\le 2\beta_0$.  Then $1.008<20/(2\beta_0)\le n/(2\beta_0)\le 2^k$,
so $k\ge 1$.  Therefore the argument above implies that for $n\ge 20$ the
bound $D\ge 0.7525\sqrt{n}$ holds.  For this lower bound to be
nontrivial we also require $0.7525\sqrt{n}\ge \log_2{n}-1$. However
this holds for all $n\ge1$ and so the lower bound on $D$ holds and is
nontrivial for all $n\ge20$.
\end{proof}

\begin{rem} A similar argument shows that there exists $\varepsilon>0$ such that if
$A \subseteq \mathbb{R}^n$ is approximately convex and satisfies
$\mathcal{H}(A,\Co(A)) \ge \log_2n - \varepsilon$, then
$\diam(A) \ge 1.16 \sqrt n$ for infinitely many $n$.
\end{rem}

Finally, we improve the upper estimate for the diameter provided by 
Corollary~\ref{cor: l_pcase}.

\begin{thm} \label{thm: Eucset}
 Let $(e_i)_{i=0}^{n}$ be the unit vector basis of $\ell_2^{n+1}$.
Then, for $n \ge 4$ and $M= \sqrt{(2/\ln 2)n\log_2n}$, the set
$$A = \left\{M\sum_{t=1}^n t_i e_i + E_{n-1}(t_1,\dots,t_n)e_0: (t_1,\dots,t_n)
\in \Delta_{n-1}\right\}$$
is approximately convex and satisfies the following:
$$ \mathcal{H}(A,\Co(A))=\log_2n \quad \text{and} \quad \operatorname{diam}(A)
\le \frac{2}{\sqrt{\ln 2}}\sqrt{n\log_2n} + \log_2n.$$
\end{thm}
\begin{rem} Theorem~\ref{thm: Eucset} is a significant improvement on
Corollary~\ref{cor: l_pcase} as it eliminates the dependence on $\varepsilon$
and reduces the exponent of $\log n$ in the estimate
for $\operatorname{diam}(A)$. When $n+1=2^k$, the set $A$ is
very nearly extremal, since in this case
 $\mathcal{H}(A,\Co(A)) \le \log_2(n+1) = C(\mathbb{R}^{n+1})$
by Theorem~\ref{thm: power2}.
\end{rem}

The proof of this result is a consequence of the solution to
a constrained optimization problem. Consider the following functional:
\begin{equation*}
I(y) = M^2 \int_0^ny(x)^2\,dx + \left(\int_0^n \phi(y(x))\,dx\right)^2,
\end{equation*} where $y(x)$ is a non-negative function defined on the
open interval $(0,n)$. (Recall that $\phi(t)=t\log_2(1/t)$.) The problem
is to minimize $I(y)$ subject to the following constraints on $y$:
$$0 \le y \le 1\quad \text{and}\quad \int_0^n y(x)\,dx=1.$$
We prove in Lemma~\ref{lem: y0} below that, for $M^2 = (2/\ln 2)n\log_2n$, $I(y)$
is minimized by $y_0 =(1/n) \chi_{(0,n)}$.

Assuming this result, let us complete the proof of Theorem~\ref{thm: Eucset}.
\begin{proof}[Proof of Theorem~\ref{thm: Eucset}] 
Clearly, $$\mathcal{H}(A,\Co(A)) \le \max_{t \in \Delta_{n-1}}E_{n-1}(t) = \log_2n.$$
To establish the reverse inequality, we show that $d(x_0,A)=\log_2n$
for $x_0 =(M/n) \sum_{i=1}^n e_i$. Observe that 
$$d(x_0,A)^2 = \min\left\{M^2\sum_{i=1}^n\left(t_i-\frac{1}{n}\right)^2 + 
E_{n-1}(t_1,\dots,t_n)^2:(t_1,\dots,t_n)
\in \Delta_{n-1}\right\},$$
and also that
$$M^2\sum_{i=1}^n\left(t_i-\frac{1}{n}\right)^2 + E_{n-1}(t_1,\dots,t_n)^2
= g(t_1,\dots,t_n)- \frac{M^2}{n},$$
where
$$g(t_1,\dots,t_n) = M^2\sum_{i=1}^n t_i^2 + \left(\sum_{i=1}^n \phi(t_i)\right)^2.$$
Hence
$$d(x_0,A)^2 = \min\{g(t_1,\dots,t_n): (t_1,\dots,t_n) \in \Delta_{n-1}\}- \frac{M^2}{n}.
$$
But $g(t_1,\dots,t_n) = I(\tilde g)$, where
$\tilde g(x) = \sum_{k=1}^n t_k \chi_{[k-1,k)}.$
(Note that $\tilde g(x)$ satisfies the constraints for the optimization problem.)
Since $I(y)$ is minimized by $y_0= (1/n)\chi_{(0,n)}$
(see Lemma~\ref{lem: y0}), we get
$$g(t_1,\dots,t_n) = I(\tilde g) \ge I(y_0) = g(1/n,\dots,1/n).$$
Hence
$$d(x_0,A)^2 = g(1/n,\dots,1/n) - \frac{M^2}{n} = (\log_2n)^2.$$ 
Thus, $\mathcal{H}(A,\Co(A)) \ge \log_2n$. The estimate for $\operatorname{diam}(A)$
is straightforward.
\end{proof} 

The next four lemmas solve the constrained optimization problem.

\begin{lem} \label{lem: existencey0} Let $M>0$ and $n \ge 1$.
There exists a 
 right-continuous non-increasing function $y_0$ on $(0,n)$ which 
solves the constrained optimization problem. \end{lem}
\begin{proof} Let $m$ be the infimum of $I(y)$ taken over all
$y$ which satisfy the constraints. There exist $y_n$ ($n \ge 1$) satisfying
the constraints such that $I(y_n) \rightarrow m$ as $n \rightarrow \infty$.
By replacing each $y_n$ by its \textit{non-increasing rearrangement},
we may assume that each $y_n$ is right-continuous and non-increasing.
By Helly's selection theorem
(see e.g. \cite[p.\ 221]{Na-1}), we may also assume (by passing to a subsequence)
 that $y_n(x) \rightarrow \tilde y_0(x)$ pointwise. Since $0 \le y_n \le 1$,
it follows from the Bounded Convergence Theorem that $\tilde y_0$
satisfies the constraints and that $I(\tilde y_0) = \lim_n I(y_n) = m$.
Finally, let $y_0$ be the right-continuous modification of $\tilde y_0$. 
\end{proof}
\begin{lem} \label{lem: y0values} There exists $\alpha \in (0,1)$ such
that the set of values taken by $y_0$ is a subset of $\{0,1,\alpha\}$.
\end{lem}
\begin{proof}In the notation of Lemma~\ref{lem: existencey0}, we may assume
that $y_k$ is a step function minimizing $I(y)$ over all
step functions of the
 form \linebreak $\sum_{j=1}^k a_j \chi_{[(j-1)n/k,jn/k)}$  satisfying the constraints.
 A  value $\lambda \in (0,1)$ taken by $y_k$ must satisfy the 
following
Lagrange multiplier equation for a local minimum:
\begin{equation} \label{eq: Lagrangemult1} 
2M^2\lambda+\left(2\int_0^n \phi(y_k)\,dx\right)\phi'(\lambda) = 2A, \end{equation}
where $A$ is a constant. It is easily seen that this equation has at most two roots in $(0,1)$.
By the pointwise convergence of $y_k$ to $\tilde y_0$,
it follows that $y_0$ takes at most two values in $(0,1)$.
Therefore we may apply the method of Lagrange multipliers again to deduce that
these values must also satisfy \eqref{eq: Lagrangemult1} (with $y_k$ replaced by $y_0$). 
Equivalently, setting $B=\int_0^n \phi(y_0)\,dx>0$,
\begin{equation} \label{eq: Lagrangemult2}
M^2\lambda + B(\log_2(1/\lambda)-(1/\ln 2))=A. \end{equation} 

Suppose that there are two distinct roots, $\alpha$
and $\beta$, with
$0 < \alpha < \beta <1$, and
suppose that $y_0$ takes one of these values, $\alpha$ say, on an interval $J$. 
(The argument is similar if $y_0$ takes the value $\beta$.) 
Let $g$ take the value $0$ on the complement of $J$, and
the values $1$ and $-1$ on the left-hand and right-hand halves of $J$,
respectively. Since $\alpha \in (0,1)$, it follows
that  $y_0 + \varepsilon g$ satisfies the constraints, provided $\varepsilon>0$
is sufficiently small. Moreover,
$$I(y_0+\varepsilon g)= I(y_0) + |J|\left(M^2- \frac{B}{(\ln 2)\alpha}\right)\varepsilon^2
+ o(\varepsilon^2).$$
Since $y_0$ minimizes $I$,
 \begin{equation} \label{eq: alphacon}
M^2- \frac{B}{(\ln 2) \alpha} \ge 0. \end{equation}
To derive a contradiction, suppose that $y_0$
also takes the value $\beta$ on an 
interval. Then, by the same argument,
\begin{equation*} 
M^2- \frac{B}{(\ln 2) \beta} \ge 0. \end{equation*}
Since \eqref{eq: Lagrangemult2} is satisfied by $\lambda= \alpha$ and $\lambda=\beta$,
the Mean Value Theorem implies the existence of $\gamma \in (\alpha,\beta)$ such that
\begin{equation*}
M^2- \frac{B}{(\ln 2) \gamma} = 0. \end{equation*}
Thus, \begin{equation*}
M^2- \frac{B}{(\ln 2) \alpha} < M^2- \frac{B}{(\ln 2) \gamma} = 0. \end{equation*}
But this contradicts \eqref{eq: alphacon}. Thus, $y_0$ cannot take the value
$\beta$, which completes the proof.
\end{proof}
\begin{lem} \label{lem: y0not1} Suppose that $n\ge4$ and that
$$5n < M^2\le \frac{2}{\ln 2}n\log_2n.$$  
Then $y_0$ does not take
the value $1$. \end{lem}
\begin{proof} For $n\ge4$, we have \begin{equation} \label{eq: I(y0)bound}
I(y_0) \le I\left(\frac{1}{n}\chi_{[0,n]}\right) = \frac{M^2}{n} + (\log_2 n)^2 \le \frac{2}{\ln 2}\log_2n + 
(\log_2n)^2 <\frac{M^2}{2}.
\end{equation}
Suppose that $y_0$ takes the value $1$ on $[0,x]$ and the nonzero 
value
$k \in (0,1)$ on an interval of length $(1-x)/k \le n-x$. If $k\ge1/2$ then
$I(y_0) \ge (1/2)M^2$, which contradicts \eqref{eq: I(y0)bound}.
So we may assume that $k \in (0,1/2)$. Now
$$ I(y_0) = M^2(x+k(1-x)) +((1-x)\log_2(1/k))^2.$$
So \begin{align*}
\frac{\partial I(y_0)}{\partial x}
&=M^2(1-k) - 2\log_2(1/k)^2(1-x)\\ 
&\ge \frac{M^2}{2} - 2\log_2(1/k)^2k(n-x)\\
\intertext{(since $(1-x) \le k(n-x)$)}
&\ge \frac{M^2}{2} - 2 \left(\max_{0\le k \le 1/2} k\log_2(1/k)^2\right)n\\
&\ge\left(\frac{5}{2}- \frac{8}{e^2(\ln 2)^2}\right)n >0.
\end{align*}
Since $I(y_0)$ minimizes $I(y)$, it follows that $x=0$, as desired.
\end{proof}
\begin{lem} \label{lem: y0} Suppose that $n \ge 4$ and that
$M^2=(2/\ln 2)n\log_2n$. Then $y_0 = (1/n) \chi_{(0,n)}$ and
$I(y_0) = 2\log_2n + (\log_2n)^2$. \end{lem}
\begin{proof} By Lemma~\ref{lem: y0not1}, $y_0$ takes only one nonzero value
$k \in [1/n,1)$ on an interval of length $1/k$. So
$I(y_0)= M^2k + (\log_2(1/k))^2.$ Thus,
\begin{align*}
\frac{\partial I(y_0)}{\partial k}&=M^2 - \frac{2(\log_2(1/k))}{(\ln2)k}\\
&= \frac{2}{\ln 2}\left(n\log_2n-\frac{1}{k}\log_2(1/k)\right) \ge 0, \end{align*}
with equality if and only if $k=1/n$. Since $y_0$ minimizes $I(y)$,
it follows that $k=1/n$, which gives the
result.
\end{proof}
\begin{rem} Setting $M^2=6n$ in Lemma~\ref{lem: y0} yields an approximately convex set
set $A \subset \mathbb{R}^{n+1}$ with $\operatorname{diam}(A)=O(\sqrt n)$ and 
$\mathcal{H}(A,\Co(A)) \ge \log_2n-c\log_2\log_2n$ for some constant $c$. \end{rem}

\section{Lower bounds in  spaces of type $p$} \label{sec: diamlowertype}

First we recall the notion of \textit{type}. 
In the following
definition $(\varepsilon_i)_{i=1}^\infty$ is a sequence of
independent Bernoulli random variables, with $P(\varepsilon_i=1)
=P(\varepsilon_i=-1)=1/2$, defined on a probability space $(\Omega,\Sigma,P)$.
The \textit{expected value} of a random variable $Y$ is denoted $\mathbb{E}Y$.
\begin{definition} Let $1 \le p \le 2$. A normed space $X$ is of
type $p$ 
if there exists a constant $T_p(X)$ (the `type $p$ constant') such that
\begin{equation*}
\left(\mathbb{E} \left\|\sum_{i=1}^n \varepsilon_i x_i\right\|^p\right)^{1/p}
\le T_p(X)\left(\sum_{i=1}^n \|x_i\|^p\right)^{1/p}\end{equation*}
for all $n\ge1$ and for all choices of $x_i \in X$ ($1 \le i \le n$).
\end{definition}
The following theorem  can be deduced from (and in 
fact is essentially equivalent to) \cite[Thm.\ 3.6]{Ca-Pa-1}.
For completeness we give a short direct proof.
We show in
Corollary~\ref{cor: l_pcasesharp} below that the exponent
of $(p-1)/p$ in this theorem is sharp.
\begin{thm} \label{thm: typep}
Let $1<p \le 2$ and let $X$ be a normed space of type $p$. Suppose that
$A \subseteq X$ is approximately convex. Let $D = \operatorname{diam}(A)$
and let $d = \mathcal{H}(A,\Co(A))$. Then, provided $d \ge 2$, we have
 \begin{equation} \label{eq: Dlowerbound}
D \ge \frac{8^{1/p}}{16T_p(X)} (2^d)^{(p-1)/p}
\end{equation}
\end{thm}
\begin{proof} 
 We may assume (cf.\ Theorem~\ref{th: lowbound3}) that
 $0 \in \Co(A)$, that 
$d=d(0,A)$, and that
$\|a\| \le D$ for all $a \in A$. Since $0 \in \Co(A)$ there exist
$m \ge 1$,
$a_i \in A$ and  $p_i >0$ ($1 \le i \le m$), with $\sum_{i=1}^{m} p_i =1$
and
$\sum_{i=1}^m p_ia_i=0$.

Let $(Y_j)_{j=1}^\infty$ be a sequence of independent
identically distributed $X$-valued
random variables defined by
\begin{equation*}
P(Y_j=a_i)=p_i\qquad (1 \le i \le m).
\end{equation*} Then $\|Y_j(\omega)\|\le D$ 
($\omega \in \Omega$) and 
$\mathbb{E}$ $Y_j = \sum_{i=1}^m p_ia_i=0$. Thus, 
\cite[Prop. 9.11]{Le-Ta-1} yields (for each $n$)
\begin{align*}
\left(\mathbb{E}\left\|\sum_{i=1}^n Y_i
\right\|^p\right)^{1/p} &\le 2T_p(X) \left(\sum_{i=1}^n \mathbb{E}
\|Y_i\|^p\right)^{1/p}\\
&\le 2T_p(X)n^{1/p}D.
\end{align*}
So there exist $b^n_i \in A$ ($1 \le i \le n$) with
\begin{equation}\label{eq: CAP}
\left\|\frac{1}{n}\sum_{i=1}^n b^n_i\right\| \le 2T_p(X) n^{(1-p)/p}D
\end{equation}
Since $A$ is approximately convex,
$$d\left(\frac{1}{n}\sum_{i=1}^n b^n_i,A\right)\le \kappa(n-1) \le \log_2n +1.$$ So
\begin{align*}
d(0,A) &\le \left\|\frac{1}{n}\sum_{i=1}^n b^n_i\right\| + 
d\left(\frac{1}{n}\sum_{i=1}^n b^n_i,A\right)\\
&\le 2T_p(X) n^{(1-p)/p}D + \log_2n +1
\end{align*}
Put $n=2^{[d]-2}$ (noting that $n \ge 1$ since $d \ge 2$ by assumption)
 so that $\log_2n+1 \le d-1$. Then
\begin{equation*}
d=d(0,A)\le 2T_p(X)D (2^{d-3})^{(1-p)/p} + d-1, \end{equation*}
which yields \eqref{eq: Dlowerbound}. \end{proof}
\begin{rem}\label{rem: CAP}  \eqref{eq: CAP} and its probabilistic proof
are from \cite{Br-1}. It is proved in \cite{Br-1} that $X$ has the
\textit{convex approximation property} if and only if 
$X$ has type $p$ for some $p>1$. When $X$ is a Hilbert space,
 Theorem~\ref{th: bestconv} above gave a deterministic proof 
of \eqref{eq: CAP} with the sharp constants.
\end{rem} 
\begin{cor} Let $X$ be a Banach space. The following are equivalent:
\begin{itemize}\item[(a)] $X$ is B-convex;
\item[(b)] there exists $c>0$ such that for every approximately convex set $A\subseteq X$, 
we have \begin{equation} \label{eq: explowerbound}
\operatorname{diam}(A) \ge c \exp(c\mathcal{H}(A,\Co(A))). \end{equation} \end{itemize}
\end{cor}
\begin{proof} It is known that $X$ is B-convex if and only if $X$ has type
$p$ for some $p>1$ \cite{Pi-1}. Thus,  (a)$\Rightarrow$(b) follows
from Theorem~\ref{thm: typep}. Now suppose that $X$ is not B-convex. By definition
(see Section~\ref{sec: appconvexsets}),
$X$ contains `almost isometric' copies of $\ell_1^n$ for all $n$. So by
Remark~\ref{rem: l1case} $X$ contains approximately convex sets $A_n$
such that $\mathcal{H}(A_n,\Co(A_n)) \ge \log_2n -1$ and 
$\operatorname{diam}(A_n) \le C\log_2n$, where $C$ is an absolute constant.
Clearly, \eqref{eq: explowerbound} cannot hold in $X$, and  so (b)$\Rightarrow$(a).
\end{proof}
\begin{rem} The above result is essentially equivalent to \cite[Thm.\ 3.7]{Ca-Pa-1},
which was first obtained in \cite{Lar-1}. \end{rem}
The following corollary is a partial converse to Corollary~\ref{cor: l_pcase}.
When combined with the latter it shows
that the factor $n^{(p-1)/p}$ in Corollary~\ref{cor: l_pcase} 
and the exponent of $(p-1)/p$ in Theorem~\ref{thm: typep} are both sharp.
\begin{cor} \label{cor: l_pcasesharp}
 Let $1 < p < \infty$. There exists a constant $c_p>0$ such that
if $A \subseteq L_p(0,1)$ is approximately convex 
and satisfies \linebreak $\mathcal{H}(A,\Co(A)) \ge \log_2n-1$, then 
\begin{equation*}
\diam(A)\ge \begin{cases}  c_p n^{(p-1)/p} &\text{($1< p \le 2$),}\\ 
c_p n^{1/2} &\text{($2 \le p <\infty$}). \end{cases} \end{equation*}
\end{cor}\begin{proof} It is known that $L_p(0,1)$ has type $\min(p,2)$. Setting
$d=\log_2n-1$ in Theorem~\ref{thm: typep} gives the result. \end{proof}
\section{Sets with $\diam(A)=\mathcal{H}(A,\Co(A))$} \label{sec: infinitedim}
In this section we show that there exists 
an infinite-dimensional Banach space $Y$ such that
for every prescribed diameter $D$
there exists
 an approximately convex set $A\subseteq Y$
such that 
$\operatorname{diam}(A)=\mathcal{H}(A,\Co(A))=D$.
This is clearly  ``worst possible''. 
More precisely, we shall prove the following theorem. 
\begin{thm} \label{thm: worst1} Let $M > 0$. There exist a Banach space 
$(X,\|\cdot\|)$ that is linearly isomorphic to $\ell_1$
 and an approximately convex set $A \subseteq B_M(X)$ 
such that $\mathcal{H}(A, \Co(A))= \operatorname{diam}(A)= 2M$.
\end{thm}
First observe that Theorem~\ref{thm: worst1} admits the following reformulation in terms of
$\varepsilon$-convex sets.
\begin{thm} \label{thm: worst2} Let $\varepsilon>0$.
There exist a Banach space 
$(X,\|\cdot\|)$ that is linearly isomorphic to $\ell_1$
 and an $\varepsilon$-convex set $A' \subseteq B(X)$ 
such that $\mathcal{H}(A', \Co(A'))= \operatorname{diam}(A')= 2$.
\end{thm} \begin{proof} Let $M=1/\varepsilon$ and let $X$ and $A$ satisfy
the conclusion of
Theorem~\ref{thm: worst1}. Then $A'= \varepsilon A$ has the required properties.
\end{proof}

The following lemma is known \cite{Ca-Pa-1}, but
for completeness
we outline the proof.
\begin{lem} \label{lem: appjensen}
 Suppose that $A \subseteq X$ is approximately Jensen-convex.
Then $A$ is $2$-convex. In particular,
$(1/2)A$ is approximately convex. \end{lem}
\begin{proof} Let $f(x) = d(x,A)$ ($x \in X$). Then $f$ is a continuous
approximately Jensen-convex function, i.e.\
\begin{equation*}
f\left(\frac{x+y}{2}\right) \le \frac{1}{2}(f(x)+f(y))+1.
\end{equation*}
By \cite{Ng-Ni-1} $f$ is a $2$-convex function, which implies that $A$
is a $2$-convex set. \end{proof}

Lemma~\ref{lem: appjensen} shows that Theorem~\ref{thm: worst1} is equivalent to
the following result. 
\begin{thm} \label{thm: worst3}
Let $M \in \mathbb{N}$. There exist a Banach space 
$(X,\|\cdot\|)$ that is linearly isomorphic to $\ell_1$
 and an approximately Jensen-convex set $A \subseteq B_M(X)$ 
such that $\mathcal{H}(A, \Co(A))= \operatorname{diam}(A)= 2M$.
\end{thm}

\begin{rem} The restriction $M\in \mathbb{N}$ is made only to simplify notation in the proof. 
 Clearly
the result will hold for all $M>0$
by scaling. \end{rem}

The rest of the paper is devoted to the 
lengthy proof of Theorem~\ref{thm: worst3}.
To construct the
 space $X$ appearing in the conclusion of the theorem,
let us begin with the `tree-like' combinatorial 
structure which will form a Schauder basis for $X$. 
Fix $M \in \mathbb{N}$. Let $L_1=\mathbb{N}$,
and for
$n>1$
define $L_n$ recursively as follows:
$$
L_n=\{(a,b):a \in L_i,b\in L_j, i+j=n, 1\le i,j<n\}.
$$
Let $L=\cup_{n=1}^\infty L_n$, and, for $a \in L$, let
$e_a$ denote the indicator function of $\{a\}$.
For $R \subseteq L$,
let $c_{00}(R)$ denote
the vector subspace of $\ell_\infty(R)$ spanned by
the set $\{e_a: a \in R\}$.
For $x = \sum_{a\in L}\lambda_a e_a \in c_{00}$,
let $\operatorname{supp}(a) = \{a \in L: \lambda_a \ne 0\}$
 
We introduce two norms,
$\|\cdot\|_1$ and $\|\cdot\|_1'$, on $c_{00}(L)$:
$$\left\|\sum_{a \in L}\lambda_a e_a\right\|_1=\sum_{a \in L}|\lambda_a|$$
and
$$\left\|\sum_{a \in L}\lambda_a e_a\right\|_1'
=M\sum_{a \in L_1}|\lambda_a|+\sum_{a \in L\setminus L_1}
|\lambda_a|.$$
Note that $\|\cdot\|_1$ is the usual $\ell_1$ norm
and that $\|\cdot\|_1'$ is a weighted $\ell_1$ norm
with respect to the basis $\{e_a:a \in L\}$.
A linear mapping  $T:c_{00}(L) \rightarrow c_{00}(L)$ 
is defined (extending linearly) thus:
\begin{equation*}
T(e_a)=\begin{cases} 0 & \text{if $a \in L_1$},\\
\dfrac{e_b+e_c}{2} & \text{if $a \in \cup_{n=2}^\infty L_n$
and $a=(b,c)$}. \end{cases} \end{equation*}
Note that $T(c_{00}(L_n))\subseteq c_{00}(\cup_{k=1}^{n-1}L_k)$
and that $T^n(x)=0$ for all $x \in c_{00}(L_n)$. Hence
$S=I-T$ is an invertible operator on $c_{00}(L)$
with inverse $S^{-1}=\sum_{k=0}^\infty T^k$. Note also that
$$\sum_{a\in L}T(x)(a) \le  \sum_{a \in L} x(a)$$
if $x(a) \ge 0$ for all $a \in A$, with equality if $\operatorname{supp}(x)
\subset \cup_{n=2}^\infty L_n$.

Define a norm $\|\cdot\|$ on $c_{00}(L)$ thus:
$$\|x\|=\inf\{M\|y\|_1 + \|S^{-1}(z)\|_1':x=y+z\}\qquad(x\in c_{00}(L)).
$$
Let $(X,\|\cdot\|)$ be the completion of $(c_{00}(L),\|\cdot\|)$
and let $A=\{e_a:a\in L\}\subseteq X$. 

The verification that
$X$ and $A$ satisfy the conclusion of Theorem~\ref{thm: worst3} will
be broken down into four lemmas.

\begin{lem} \label{lem: Fandphi} Suppose that $F \in B(X^*)$. Then 
the mapping
 $\phi:L\rightarrow \mathbb{R}$
defined by $\phi(a) = F(e_a)$ satisfies the following:
\begin{itemize}
\item[(a)] $|\phi(a)| \le M$ for all $a \in L$;
\item[(b)]$$ \left|\phi(a)-\frac{\phi(b)+\phi(c)}{2}\right| \le 1$$
for all $a=(b,c) \in \cup_{n=2}^\infty L_n$. \end{itemize}
Conversely, every $\phi$ which satisfies (a) and (b) corresponds to a
unique $F \in B(X^*)$.
\end{lem}
\begin{proof}
>From the definition of $\|\cdot\|$ we see that $F \in B(X^*)$
if and only if
\begin{equation} \label{eq: Funit}
|F(x)| \le \min(M\|x\|_1,\|S^{-1}(x)\|_1') \qquad(x \in c_{00}(L)).
\end{equation}
Indeed, if $F$ satisfies \eqref{eq: Funit}, then for every $x \in c_{00}(L)$, we have
\begin{align*} \|x\|&= \inf \{M\|y\|_1 + \|S^{-1}(z)\|_1':x=y+z\}\\
&\ge\inf\{F(y) + F(z):x=y+z\} = F(x), \end{align*}
and so $\|F\|\le 1$. Conversely, if $\|F\| \le 1$, then
$$ F(x) \le \|x\| \le \min(M\|x\|_1,\|S^{-1}(x)\|_1'),$$
and \eqref{eq: Funit} is satisfied.
  
The condition $|F(x)| \le M\|x\|_1$ is  clearly equivalent to (a). 
Since $\|\cdot\|_1$ is
a weighted $\ell_1$ norm, the condition
$|F(x)| \le \|S^{-1}(x)\|_1'$ is equivalent to
the condition \begin{equation} \label{eq: F(x)}
|F(S(e_a))| \le \|e_a\|_1' \qquad(a \in L). \end{equation} 
Suppose that
$a \in L_1$. Then $S(e_a)=e_a$ and $\|e_a\|_1'=M$, and so
\eqref{eq: F(x)} becomes $|\phi(a)|\le M$. Now suppose that
$a=(b,c) \in \cup_{n=2}^\infty L_n$. Then $S(e_a)=e_a-(1/2)(e_b+e_c)$
and $\|e_a\|_1'=1$, and so \eqref{eq: F(x)} becomes
$$\left|\phi(a)-\frac{\phi(b)+\phi(c)}{2}\right| \le 1,$$
and so (b) is satisfied.
Conversely, if $\phi$ satisfies (a) and (b), then
the mapping $F(e_a)=\phi(a)$ will extend linearly 	to an element of $B(X^*)$.
\end{proof}
\begin{rem}\label{rem: ell1} From the description of $X^*$ it follows that
$$\frac{1}{2}\|x\|_1\le \|x\| \le M \|x\|_1.$$ So $(X,\|\cdot\|)$
is isomorphic to $\ell_1$ (and the Banach-Mazur distance from
$X$ to $\ell_1$ is at most $2M$). \end{rem}
\begin{lem}\label{lem: extendphi} Suppose that $E \subseteq L$ has the property 
that whenever
 $a=(b,c) \in E$, then $b,c \in E$. If $\phi_0:E \rightarrow
[-M,M]$ satisfies
\begin{equation}\label{eq: phi_0}
\left|\phi_0(a)-\frac{\phi_0(b)+\phi_0(c)}{2}\right| \le 1
\end{equation}
for all $a=(b,c) \in E$, then $\phi_0$ admits an extension
$\phi:L \rightarrow [-M,M]$ satisfying \begin{equation}\label{eq: phi}
\left|\phi(a)-\frac{\phi(b)+\phi(c)}{2}\right| \le 1
\end{equation}
for all $a=(b,c) \in \cup_{n=2}^\infty L_n$. \end{lem}
\begin{proof} We define $\phi$ recursively. First define
$\phi$ from $L_1$ into $[-M,M]$ to be an
arbitrary extension of the restriction of $\phi_0$ to $L_1$.
Suppose that $n>1$ and that $\phi$ has been defined on
$\cup_{k=1}^{n-1}L_k$ to extend the restriction of
 $\phi_0$ to  $\cup_{k=1}^{n-1}L_k$. Let $a=(b,c) \in L_n$.
Then $b,c \in \cup_{k=1}^{n-1}L_k$, and so $\phi(b)$
and $\phi(c)$ have already been defined.
If $a \in E$, then $b,c \in E$,
and so $\phi(b)=\phi_0(b)$
and $\phi(c)=\phi_0(c)$.
It follows from \eqref{eq: phi_0} that \eqref{eq: phi} will be satisfied
with $\phi(a)=\phi_0(a)$. If $a \notin E$, define
$\phi(a)= (1/2)(\phi(b)+\phi(c))$, so that \eqref{eq: phi} is
trivially satisfied.
This completes the definition of $\phi$ on $L_n$.
\end{proof}
Now fix $a \in L$ and let $E_a=\cup_{n=0}^\infty 
\operatorname{supp}(T^n(e_a))$ ($=\cup_{n=0}^{N-1} 
\operatorname{supp}(T^n(e_a))$ for $a \in L_N$).
For $d \in E_a$, we define the {\em $a$-order} of
$d$, denoted $\mathbf{o}_a(d)$, thus:
$$\mathbf{o}_a(d)=\min\{n\ge0:d\in\operatorname{supp}(T^n(e_a))\}.$$
\begin{lem}\label{lem: phi(d)} Given $a \in L$, there exists $\phi: L
\rightarrow [-M,M]$ satisfying \eqref{eq: phi} such that \begin{equation}
\label{eq: phi(d)-M}\phi(d)=-M \qquad(d \in L_1\setminus E_a),\end{equation}
\begin{equation}
\label{eq: phi(d)L1} \phi(d) = \max(M-\mathbf{o}_a(d),-M)\qquad (d \in L_1 \cap E_a),
 \end{equation}
and
\begin{equation}\label{eq: phi(d)}
\phi(d) \ge \max(M-\mathbf{o}_a(d),-M)\qquad (d \in E_a).
\end{equation}
\end{lem}
\begin{proof} First we define  a mapping
 $\phi_0:E_a \cup L_1\rightarrow [-M,M]$. 
If $d \in L_1\setminus E_a$, let $\phi_0(d)=-M$, and if
$d \in E_a\cap L_1$, let
$$\phi_0(d)=\max(M-\mathbf{o}_a(d),-M),$$
so that \eqref{eq: phi(d)-M} and \eqref{eq: phi(d)L1} are satisfied.
Now extend to the rest of $E_a$ recursively as follows.
Suppose that $n >1$ and that
 $\phi_0$ has been defined on $E_a \cap (\cup_{k=0}
^{n-1}L_k)$ to satisfy \eqref{eq: phi} and \eqref{eq: phi(d)}.
 Let $d\in E_a \cap L_n$. Then $d=(b,c)$
for some $b,c \in E_a \cap (\cup_{k=0}
^{n-1}L_k)$. Define
$$ \phi_0(d)= \min\left(M,\frac{\phi_0(b)+\phi_0(c)}{2}+1\right).$$
If $\phi_0(d) = M$, then, as $\phi_0(b)\le M$ and $\phi_0(c) \le M$, we  have
$$\phi_0(d) = M \le \frac{\phi_0(b)+\phi_0(c)}{2} + 1 \le \frac{M+M}{2} +1 =M+1,$$
so that 
$$
\left|\phi_0(d)-\frac{\phi_0(b)+\phi_0(c)}{2}\right| \le 1,
$$
i.e., \eqref{eq: phi} is satisfied by $d=(b,c)$. Also, if $\phi_0(d)=M$, then
\eqref{eq: phi(d)} is trivially satisfied.

On the other hand, if
$\phi_0(d)=(\phi_0(b)+\phi_0(c))/2+1$, then \eqref{eq: phi} is trivially satisfied
by $d=(b,c)$. In order to verify \eqref{eq: phi(d)}, suppose that
$\mathbf{o}_a(d)=k$. Then both $\mathbf{o}_a(b) \le k+1$
and $\mathbf{o}_a(c) \le k+1$.  Moreover, both $b$ and $c$ satisfy
\eqref{eq: phi(d)} by the recursive hypothesis.
Thus,
\begin{align*}
\phi_0(d) &= \frac{\phi_0(b)+\phi_0(c)}{2}+1\\
&\ge \frac{\max(M-\mathbf{o}_a(b),-M)+\max(M-\mathbf{o}_a(c),-M)}{2}+1\\
&\ge \frac{\max(M-(k+1),-M)+\max(M-(k+1),-M)}{2}+1\\
&\ge \max(M-k,-M)\\
&=\max(M-\mathbf{o}_a(d),-M). \end{align*}
Thus, \eqref{eq: phi(d)} is satisfied by $d$, which completes the recursive
definition of $\phi_0$.
 Now $\phi_0$ and $E_a \cup L_1$ (replacing $E$)
satisfy the hypotheses
of Lemma~\ref{lem: extendphi}. Let $\phi$ be the extension of $\phi_0$ given by 
Lemma~\ref{lem: extendphi}.
\end{proof}
The following lemma completes the proof of Theorem~\ref{thm: worst3}.
\begin{lem} Let $A=\{e_a:a \in L\}$. Then $A$ satisfies the following:
\begin{itemize}
\item[(i)] $A \subseteq B_M(X)$; 
\item[(ii)] A is approximately Jensen-convex;
\item[(iii)] $\mathcal{H}(A, \Co(A)) = 2M$. 
\end{itemize} \end{lem}
\begin{proof} Suppose that $F \in B(X^*)$. By Lemma~\ref{lem: Fandphi},
$|F(e_a)|\le M$ for all $a \in A$, and so (i) follows from the
Hahn-Banach Theorem. Suppose that $b,c \in A$. Then $a=(b,c) \in A$,
and by Lemma~\ref{lem: Fandphi} $$|F(e_a) - F((1/2)(e_b+e_c))| \le 1,$$ which gives (ii).
To prove (iii), note that
(i) implies that $\Co(A) \subseteq B_M(X)$
(since $B_M(X)$ is convex), and hence
$$\mathcal{H}(A, \Co(A)) \le
\operatorname{diam}(B_M(X))=2M.$$
So it suffices to prove that $\mathcal{H}(A, \Co(A)) \ge 2M$.
Fix $N \ge 1$ and choose distinct elements $a_1,\dots,a_N \in L_1$.
We shall prove that 
$$d\left(\frac{1}{N}\sum_{k=1}^N e_{a_k},A\right) \ge 2M - \varepsilon_N,$$
where $\varepsilon_N \rightarrow 0$ as $N \rightarrow \infty$.
Let $a \in L$. If $d \in E_a$ and $\mathbf{o}_a(d)=k$, then 
$T^k(e_a)(d) \ge 2^{-k}$. Since $\sum_{b\in L} T^k(e_a)(b)\le1$,
it follows that $E_k=\{d \in L: \mathbf{o}_a(d)=k\}$ has cardinality
at most $2^k$. Thus
$$\left| \bigcup_{k=0}^{2M-1} E_k\right| \le \sum_{k=0}^{2M-1} 2^k = 2^{2M}-1.$$
Let $\phi: L \rightarrow [-M,M]$ be the function
associated to $a$  
 defined  in Lemma~\ref{lem: phi(d)}, and
let $F\in B(X^*)$
be the linear functional corresponding to $\phi$.
If $a_i \in L_1 \setminus E_a$, then $\phi(a_i)=-M$ by \eqref{eq: phi(d)-M}.
If $a_i \in E_a$ and $\mathbf{o}_a(a_i) \ge 2M$, then
$\phi(a_i)= -M$ by \eqref{eq: phi(d)L1}. Hence
if $a_i \notin G=\cup_{k=0}^{2M-1} E_k$, then $\phi(a_i)=-M$.
Moreover, $\phi(a)=M$  by \eqref{eq: phi(d)}, since $\mathbf{o}_a(a)=0$.  So
\begin{align*}
F\left(e_a - \frac{1}{N}\left(\sum_{k=1}^N e_{a_k}\right)\right) &=
\phi(a) - \frac{1}{N}\sum_{k=1}^N \phi(a_i)\\
&\ge M - \frac{1}{N}((N - |G|)(-M) + |G|M)\\
&= 2M - \frac{2}{N}|G|M\\
&\ge 2M - \frac{2^{2M+1}M}{N}, \end{align*}
and so 
$$\left\|\frac{1}{N}\sum_{k=1}^N e_{a_k}-e_a\right\| \ge 2M - \varepsilon_N,$$
where $\varepsilon_N = 2^{2M+1}M/N \rightarrow 0$ as $N \rightarrow
\infty$ as desired.
\end{proof}
\begin{thm} \label{thm: spaceY} There exists a Banach space $Y$ such that for 
every $\varepsilon>0$
there exists an $\varepsilon$-convex set $A \subseteq B(Y)$ with
$\mathcal{H}(A,\Co(A))=2$. \end{thm}
\begin{proof}
Let $X_n$ denote the space constructed above for $M=n$.
Then the $\ell_2$-sum
 $Y=(\sum_{n=1}^\infty \oplus X_n)_2$ has the
required property. \end{proof}
\begin{rem} Since $X_n$ is isomorphic to $\ell_1$ (Remark~\ref{rem: ell1}),
 it has  both the \textit{Radon-Nikod\'ym} property 
(see e.g.\ \cite{Di-Uh-1})
and the \textit{approximation property} (see e.g.\ \cite[p.\ 29]{Li-Tz-1}).
Hence
$Y=(\sum_{n=1}^\infty \oplus X_n)_2$ has the Radon-Nikod\'ym property \cite[p. 219]{Di-Uh-1} 
and (as is easily verified) the approximation property. \end{rem}

Since $C(0,1)$ is a universal space
for separable Banach spaces \linebreak (Mazur's theorem), it
satisfies the conclusion of Theorem~\ref{thm: spaceY}. So,
finally, let us 
reformulate
Theorem~\ref{thm: worst2} to
make good the claim made in Remark~\ref{rem: approxlip}.

\begin{cor} \label{cor: approxlip} Let $\varepsilon>0$. There exists a (non-negative)
 $\varepsilon$-convex
$1$-Lipschitz function on $B(C(0,1))$ such that 
$$ \sup\{|f(x)-g(x)|:x \in B(C(0,1))\}\ge1$$
for every convex function $g$. \end{cor}
\begin{proof} By Theorem~\ref{thm: worst2} there exists $A \subseteq B(C(0,1))$
such that $A$ is $\varepsilon$-convex and $\mathcal{H}(A,\Co(A))=2$.
Then $f(x)= d(x,A)$ has the required properties. \end{proof}

\end{document}